\documentclass{amsart}
\usepackage[T1]{fontenc}
\usepackage[cp1250]{inputenc}
\usepackage{amsrefs}

\usepackage{amsmath,amsthm,amssymb}
\usepackage{amscd}

\newenvironment{Proof}{\textbf{Proof.}}{$\qquad \blacksquare$\par}
\newenvironment{Proof of}[1]{\textbf{Proof #1.}}{$\qquad \blacksquare$\par}

\newcommand{\B}{\mathcal B}

\newcommand{\K}{\mathcal K}

\newcommand{\M}{Mes}

\newcommand{\al}{\gamma}
\DeclareMathOperator{\Int}{Int}

\renewcommand{\L}{\mathcal L}

\newcommand{\A}{\mathcal A}

\newcommand{\C}{\mathbb C}

\newcommand{\Z}{\mathbb Z}
\newcommand{\N}{\mathbb N}

\newcommand{\supp}{\textrm{supp}\,}

\newtheorem{thm}{Theorem}[section]\newtheorem{lem}[thm]{Lemma} \newtheorem{stw}[thm]{Proposition} \newtheorem{wn}[thm]{Corollary}
\theoremstyle{definition} \newtheorem{defn}[thm]{Definition}\newtheorem{prz}[thm]{Example}\newtheorem{uw}[thm]{Remark}

\begin{document}
 \title{On  transfer operators 
for $C^*$-dynamical systems}

\author{B. K.  Kwa\'sniewski} 

 \address{ Institute of Mathematics,  University  of Bialystok,\\
 ul. Akademicka 2, PL-15-267  Bialystok, Poland}
 \email{bartoszk@math.uwb.edu.pl}
\urladdr{http://math.uwb.edu.pl/~zaf/kwasniewski}

\begin{abstract}
The theme of the paper is  the question of existence and basic structure of transfer operators for  endomorphisms of a unital $C^*$-algebra. 
We establish a complete description of non-degenerate transfer operators,  
 characterize   complete transfer operators and  clarify  their  role in crossed product constructions. Also, we give necessary and sufficient conditions for existence  of transfer operators for commutative systems, and discuss their form for endomorphisms of $\B(H)$ which is relevant to
Kadison-Singer problem. 
\end{abstract}


   \thanks{This work was in part supported by Polish Ministry of Science and High Education grant number N N201 382634} 
   
\maketitle
\section*{Introduction}
In \cite{exel2} R. Exel introduced a notion of  \emph{transfer
operator} for $C^*$-dynamical systems as  a  natural
generalization of  the corresponding notion from classical
dynamics -- the Ruelle-Perron-Frobenius operator. His aim was to use  transfer operators in a
construction of crossed-products  associated with  irreversible
$C^*$-dynamical systems, and  one of the problematic
issues was the dependence of his construction  on the choice of  a
transfer operator which usually is not unique. This problem was to
some extent circumvented recently by A. Bakhtin and A. V.
Lebedev who introduced in \cite{Bakht-Leb} a notion of
\emph{complete transfer operator} which if it exists is a
necessarily unique non-degenerate transfer operator and  is a
sufficient tool to deal with the most important types of
crossed-products, see \cite{Ant-Bakht-Leb}, \cite{kwa-leb3}. Thus in this context, but also for future potential applications in noncommutative dynamics, cf. \cite{abl}, \cite{el}, it is essential to understand the structure and the relationship between the non-degenerate and complete transfer operators. So far, however, this topic    was not  thoroughly investigated  and   our objective is  a response to this  deficiency.
\par
We start with pushing to the limit a relation    touched upon  in \cite{exel2}*{Prop. 2.6} which results with  a complete description of  non-degenerate transfer operators via  conditional expectations, Theorem \ref{conditional expectations}. We improve investigations of \cite{Bakht-Leb} by showing that a complete transfer operator for an endomorphism $\alpha$ exists if and only if $\alpha$ has a unital kernel and  a hereditary range. This allow us to give an explicit definition of  complete transfer operators (Definition \ref{definition of complete transfers}), characterize them in terms of Hilbert $C^*$-modules (Proposition \ref{Hilbert bimodules}), and   clarify their predominant role in the theory of  crossed-products  by endomorphisms (Remark \ref{remark on crossed products}).\par
In Section \ref{commutative
section} we illustrate the  problem of uniqueness and existence of transfer operators including a brief survey of 
 commutative systems where transfer operators are related to the so-called averaging operators \cite{Ditor}, \cite{Pelczynski}. In particular, Theorem \ref{necessary sufficient conditions} presents necessary and sufficient conditions for existence respectively of non-zero and non-degenerate transfer operators. 
 \par
In the closing section  we  investigate $C^*$-dynamical systems on the
algebra $\B(H)$ of all bounded operators on a separable Hilbert
space $H$. We incline toward the point of view that an
endomorphism $\alpha:\B(H)\to \B(H)$ of  index $n$ may be
considered as a non-commutative $n$-to-one mapping, and thus we represent transfer operators for
$(\B(H),\alpha)$ as integrals on the space consisting of
$n$-elements. This description includes all possible transfer
operators in the case $n$ is finite, and all normal and certain singular operators  in
the  case $n=\infty$, Theorems \ref{endmorph of B(H)1}, \ref{integrating  transfers thm}. The complete description when
$n=\infty$ is challenging   and  is related to the
celebrated Kadison-Singer problem \cite{Kadison-Singer}, see page \pageref{pager ref for KS}.
 \section{Transfer operators and their characterizations}
\label{transfer}

Throughout the paper we let  $\A$ be a  $C^*$-algebra with an identity $1$ and
$\alpha:\A\to \A$ be an endomorphism  of this $C^*$-algebra, referring to the pair $(\A,\alpha)$ as to a \emph{$C^*$-dynamical system}. We
start with   recalling  basic definitions and facts concerning
transfer operators, cf. \cite{exel2}, \cite{Bakht-Leb}.
 \begin{defn}
A linear map $\L:\A\to \A$ is called a {\em transfer
operator} for  $(\A,\alpha)$ if it is continuous,
positive and such that
\begin{equation}
\L(\alpha(a)b) =a\L(b),\qquad a,b\in\A.
\label{b,,2}
\end{equation}
\end{defn}
If $\L$ is a transfer operator for $(\A,\alpha)$, by passing  to adjoints one gets the  symmetrized version of \eqref{b,,2}:
$$
\L(b\alpha(a)) =\L(b)a,\qquad a,b \in \A.
$$
This implies that the range of $\L$ is a two-sided ideal and $\L(1)$ is a positive central element in $\A$.
Clearly, every pair $(\A,\alpha)$ admits a transfer operator, namely  the zero operator
will do. In fact, as we shall see, there are $C^*$-dynamical systems  for which the  only
transfer operator  is the zero one.
In order to avoid  such degenerated cases it is natural to  impose  additional
conditions such as presented in the following  
\begin{stw}
[\cite{exel2}, Proposition 2.3]\label{Ex2.3}
Let $\L$ be a transfer operator for the pair $(\A,\alpha)$.
Then the following are equivalent:
\begin{itemize}
\smallskip
\item[i)]  $E = \alpha \circ \L$
is a conditional expectation from $\A$ onto $\alpha (\A)$,
\item[ii)] $\alpha\circ\L\circ\alpha =\alpha$,
\item[iii)] $\alpha (\L(1)) = \alpha (1)$.
\end{itemize}
\end{stw}

\begin{defn}[\cite{exel2}, Definition 2.3]  A transfer operator $\L$ is said to be  {\em non-degenerate} if the equivalent conditions of Proposition  \ref{Ex2.3} hold.
\end{defn} 
In connection with condition i) of Proposition \ref{Ex2.3} one easily sees that a  conditional expectation  from $\A$ onto $\alpha(\A)$ is uniquely determined by its  restriction to $\alpha(1)\A\alpha(1)$ which  yields a conditional expectation from $\alpha(1)\A\alpha(1)$ onto  $\alpha(\A)$.  This observation implies 
\begin{stw}\label{conditional expectations proposition}
A transfer operator $\L$ for  $(\A,\alpha)$ is non-degenerate iff  $E = \alpha \circ \L$
is a conditional expectation from $\alpha(1)\A\alpha(1)$ onto $\alpha (\A)$.
\end{stw}
  An essential part of structure and a  necessary condition for the existence of
non-degenerate transfer operators is  presented in the following
proposition which is an  improvement  of \cite{exel2}*{Prop. 2.5} and the corresponding results in \cite{Bakht-Leb}. For an ideal $I$ in $\A$ we denote by $I^{\bot}=\{a\in \A: aI=0\}$ its annihilator.
\begin{stw}\label{non-degenerate operators}
Suppose that there exists a non-degenerate transfer operator for $(\A,\alpha)$. Then  $\ker
\alpha$ is unital and hence $\A$ admits the decomposition
$$
\A = \ker\alpha \oplus (\ker\alpha)^\bot.
$$
Moreover, if $\L$ is a non-degenerate transfer operator for $(\A,\alpha)$, then 
\begin{itemize}
\item[i)] $
\L(\A)=(\ker\alpha)^\bot
$, in particular $1-\L(1)$ is the unit in $\ker\alpha$, 
\item[ii)]
 $\L: \alpha (\A) \to \L
(\A)$ is a  $^*$-isomorphism uniquely determined by $\alpha$. Namely, it is the inverse to
the $^*$-isomorphism
$\alpha: \L(\A)=(\ker\alpha)^\bot \to \alpha(\A)$.
\end{itemize}
\end{stw}
\begin{Proof} That $
\A = \ker\alpha \oplus \L(\A)
$ follows from \cite{exel2}*{Prop. 2.5}. Since however  $\alpha(1)$ is a projection, $\alpha(\L(1))=\alpha(1)$ and $\alpha:(\ker\alpha)^\bot \to \alpha(\A)$ is a  $^*$-isomorphism we get that $\L(1)=(\alpha|_{\L(\A)})^{-1}(\alpha(1))$ is a (central) projection in $\A$. Thus $1-\L(1)$ is the unit in $\ker\alpha$ and  $
\L(\A)=(\ker\alpha)^\bot
$. Using $\alpha\circ\L\circ\alpha =\alpha$ one sees that  $\L: \alpha (\A) \to \L
(\A)$ coincides with $(\alpha|_{(\ker\alpha)^\bot})^{-1}$.
\end{Proof}
Combining Propositions  \ref{conditional expectations proposition} and \ref{non-degenerate
operators} we get 
\begin{thm}\label{conditional expectations}
There exists a non-degenerate transfer operator for $(\A,\alpha)$ if and only
if $\ker\alpha$ is unital and there exists a conditional expectation $E:
\alpha(1)\A\alpha(1) \to \alpha(\A)$.
\par
More precisely, if $\ker\alpha$ is unital  we have  a one-to-one correspondence between
non-degenerate transfer operators $\L$ for $(\A,\alpha)$ and
conditional expectations $E$ from $\alpha(1)\A\alpha(1)$ onto
$\alpha(\A)$ established via the formulae
$$
E=\alpha\circ \L|_{\alpha(1)\A\alpha(1)},\qquad
\L(a)=\alpha^{-1}\big(E(\alpha(1)a\alpha(1))\big),\qquad a\in \A,
$$
where $\alpha^{-1}$ is the inverse to the isomorphism $\alpha:(\ker\alpha)^\bot \to
\alpha(\A)$.
\end{thm}
 \begin{Proof} In view of Propositions  \ref{conditional expectations proposition} and \ref{non-degenerate
operators}, it suffices to check that
$\L(a)=\alpha^{-1}(E(\alpha(1)a\delta(1)))$ is a
non-degenerate transfer operator which is straightforward.
\end{Proof}
 In view of the above the nicest
situation one may imagine  is that when 
$\alpha(\A)=\alpha(1)\A\alpha(1)$ which  holds iff  $\alpha(\A)$ is a hereditary subalgebra of $\A$, cf.  \cite{exel2}*{Prop. 4.1}. Then  the only 
conditional expectation  we may consider is  the identity, and thereby (if $\ker\alpha$ is unital) there is a unique
non-degenerate transfer operator for $\alpha$. Operator arising in this manner obviously satisfies 
$$
\alpha (\L(a)) =\alpha(1)a\alpha(1),\qquad a\in\A,
$$
and hence it coincides with the one called by I. V. Bakhtin and
A. V. Lebedev \cite{Bakht-Leb}, \cite{kwa-leb3}*{2.3},  a \emph{complete transfer operator}. Thus in contrary to \cite{Bakht-Leb}, \cite{kwa-leb3}, we may  give a definition of this notion in terms intrinsic to the system $(\A,\alpha)$.
 \begin{defn}\label{definition of complete transfers}
 We shall say that  $(\A,\alpha)$ \emph{admits a complete transfer operator} if
$\alpha:\A \to \A$ has a unital kernel and  a hereditary range. Then, by Theorem \ref{conditional expectations}, there is a unique
non-degenerate transfer operator for $(\A,\alpha)$ given by
 $$
 \L(a)=\alpha^{-1}(\alpha(1)a\alpha(1)),\qquad a\in \A,
$$
where $\alpha^{-1}$ is the inverse to the isomorphism $\alpha:(\ker\alpha)^\bot \to
\alpha(\A)$. We shall call  $\L$ 
a \emph{complete transfer operator} for $(\A,\alpha)$.
 \end{defn}

We now characterize complete transfer operators via Hilbert $C^*$-modules which will allow us to reveal the relationship between the various crossed products.
We recall \cite{fmr}*{Ex. 16}, \cite{kwakwa}, that  there is a natural structure of a  $C^*$-correspondence   over $\A$ on the space 
$
E:=\alpha(1) \A
$
  given by
$$
a \cdot x :=\alpha(a)x,\quad x\cdot a:= xa,\quad \textrm{ and}\quad \langle x,y\rangle_\A:=x^*y, \qquad x,y\in E,\,\, a\in \A.  
$$
So in particular, $E$ is a right Hilbert $\A$-module  equipped with a left action by adjointable maps, cf. \cite{lance}. Such objects are sometimes called Hilbert bimodules  \cite{fmr}. However there are reasons to reserve  this term  for an object with an additional structure. Namely, as in \cite{bms}, \cite{aee}, we adopt the following 
\begin{defn}
 A \emph{Hilbert $\A$-bimodule} is a Banach space $E$ which is both right Hilbert $\A$-module and left Hilbert $\A$-module with $\A$-valued inner products $\langle\cdot ,\cdot \rangle_\A$ and ${_\A}\langle \cdot ,\cdot \rangle$ connected via  the so-called imprimitivity condition 
$$
 x \cdot \langle y ,z \rangle_\A = {_\A\langle} x , y  \rangle \cdot z, \qquad \textrm{for all }\,\,\, \,x,y,z\in E.
$$
\end{defn}
Now, let $E$ be a $C^*$-correspondence  of a $C^*$-dynamical system $(\A,\alpha)$.
\begin{stw}\label{Hilbert bimodules}
 There exists a left $\A$-valued inner product ${_\A}\langle \cdot ,\cdot \rangle$ making $E$ (with its predefined left action) into   a Hilbert bimodule if and only if 
$(\A,\alpha)$ admits a complete transfer operator. Moreover, if $\L$ is a complete transfer operator, then 
$$
{_\A\langle} x , y  \rangle=\L(xy^*),\qquad \L(a)={_\A\langle} \alpha(1)a , \alpha(1)  \rangle\qquad  \quad x,y\in E, a\in \A.
$$
\end{stw}
\begin{Proof}
If $\L$ is a complete transfer operator, then ${_\A\langle} x , y  \rangle:=\L(xy^*)$ makes $E$ a left semi-inner $\A$-module by \eqref{b,,2} and positivity of $\L$. Since, for $x \in E=\alpha(1)\A$ we have  $xx^*\in \alpha(1)\A\alpha(1)=\alpha(\A)$ and $\L$ is faithful on $\alpha(\A)$,  ${_\A}\langle \cdot ,\cdot \rangle$ is non-degenerate.  
Furthermore, for $x,y,z\in E$
$$
{_\A\langle} x , y  \rangle \cdot z=\alpha(\L(xy^*))z=  \alpha(1)xy^*\alpha(1)z=xy^*z= x \cdot \langle y ,z \rangle_\A.
$$
Hence $E$ is a Hilbert $\A$-bimodule. Conversely, if $E$ is a Hilbert $\A$-bimodule and  $\L(a):={_\A\langle} \alpha(1)a , \alpha(1)  \rangle$, $a\in \A$, then 
$$
 \L(\alpha(b)a))={_\A\langle} \alpha(1)\alpha(b)a , \alpha(1)  \rangle={_\A\langle}  b \cdot \alpha(1)a , \alpha(1)  \rangle =b {_\A\langle}    \alpha(1)a , \alpha(1)  \rangle=b \L(a).
$$
and 
\begin{align*}
 \alpha(\L(a))&=\alpha({_\A\langle} \alpha(1)a , \alpha(1)  \rangle)=\alpha({_\A\langle} \alpha(1)a , \alpha(1)  \rangle)\alpha(1)= \\
&= {_\A\langle} \alpha(1)a , \alpha(1)  \rangle \cdot \alpha(1)=\alpha(1)a\cdot {\langle} \alpha(1) , \alpha(1)  \rangle_\A=\alpha(1)a\alpha(1).
 \end{align*}
 Moreover, using the fact  that $\A$  acts on $E$ from the  right via operators adjointable with respect to ${_\A}\langle \cdot ,\cdot \rangle$, cf. \cite{bms}*{Rem 1.9}   we get  $$
 \L(aa^*)={_\A\langle} \alpha(1)aa^* , \alpha(1)  \rangle={_\A\langle}   \alpha(1)a , \alpha(1)a  \rangle \geq 0.
 $$
 Hence $\L$ is the complete transfer operator for $(\A,\alpha)$.
\end{Proof}
\begin{uw}\label{remark on crossed products}
If $(\A,\alpha)$ admits a complete transfer operator and $E$ is the corresponding bimodule, then  similarly as in  \cite{fmr}*{Ex. 16}, cf. \cite{kwa-leb1}, \cite{kwakwa}, one may  see that there is a one-to-one correspondence between covariant representations of $(\A,\alpha)$, see  \cite[3.1]{kwa-leb3}, \cite[2.4]{Ant-Bakht-Leb} and covariant representations of the bimodule $E$, see \cite{aee}*{2.1}. Thus the \emph{crossed product} $\A\rtimes_\alpha \Z$ considered in  \cite{Ant-Bakht-Leb} (and as a particular case in \cite{kwa-leb3}) coincides with the \emph{crossed product $\A\rtimes_E \Z$ by the Hilbert bimodule} $E$ introduced in  \cite{aee}.  Furthermore,  by \cite[Thm. 4.15]{Ant-Bakht-Leb}, these algebras coincide with \emph{Exel's crossed product $\A\rtimes_{\alpha,\L}\Z$}, and by \cite[Prop. 3.7]{katsura1}, with \emph{Katsura's version of the Cuntz-Pimsner algebra} $\mathcal{O}_E$. The characteristic feature of the considered algebra is that  $\A$ embeds in it as the core - the fixed point algebra of the dual action, cf. \cite[Thm. 3.1]{aee}. This means that  every   (semisaturated) partial isometric crossed-product   may be reduced  to the construction based on a $C^*$-dynamical system  with a complete transfer operator, cf. \cite{Ant-Bakht-Leb}, \cite{kwakwa}. 
\end{uw}

 \section{Transfer operators for commutative systems} \label{commutative
section}
In  this section we assume that $\A$ is commutative. Then    $\A=C(X)$ where  $X$ is a compact
Hausdorff space and it is well known, see for instance \cite[Thm. 2.2]{kwa-leb2}, that $\alpha$ is of the form
\begin{equation}\label{endomorphisms}
\alpha (a)(x)=\left\{ \begin{array}{ll} a(\al(x))& ,\ \ x\in
\Delta\\
0 & ,\ \ x\notin \Delta \end{array}\right. ,\qquad a\in C(X),
\end{equation}
 where $\Delta$ is a clopen subset of $X$ and $\al:\Delta\to X$ is a
continuous map.  For any compact Hausdorff space $Y$ we let  $\M(Y)$
denote the space of all finite regular positive Borel measures on
$Y$ endowed with the weak$^*$ topology.
  \begin{stw}
 \label{transfer equation for commutative stw1}
  A mapping  $\L$ is a transfer operator for the commutative system
$(\A,\alpha)$ if and only if it is of the form
\begin{equation}\label{transfer equation for commutative1}
\L(a)(x)=\begin{cases}
\int_{\al^{-1}(x)} a(y) d\mu_x (y)&,  x \in
\Int(\al(\Delta))\\
0   &, x\in \overline{X\setminus \al(\Delta)}
\end{cases}
\qquad a \in \A,
\end{equation}
where the measure $\mu_x$ is supported on  $\al^{-1}(x)$ and the mapping
$$
\textrm{\emph{Int}}(\al(\Delta)) \ni x \longmapsto \mu_x \in  \M(\Delta)
$$
 is continuous and vanishing at the infinity in the sense that for every 
$a\in \A$ and every $\varepsilon >0$  the set
$
\{x:   |\int_{\al^{-1}(x)} a(y) d\mu_x (y)|\geq \varepsilon \}$  is compact in
 $\textrm{\emph{Int}}(\al(\Delta))$.
\\
In particular,
\begin{itemize}
\item[i)] $\L$ is non-degenerate iff $\al(\Delta)$ is open and $\mu_x$
is a probability measure for $x\in \al(\Delta)$,
\item[ii)] $\L$ is  complete  iff $\al(\Delta)$ is open,
$\al:\Delta\to \al(\Delta)$ is a homeomorphism, and then
$$
\L(a)(x)=\begin{cases}
a(\al^{-1}(x)) &,  x \in \al(\Delta)\\
0   &, x\in X\setminus \al(\Delta)
\end{cases}
\qquad a \in \A.
$$
\end{itemize}
 \end{stw}
\begin{Proof}
It is well known, cf. \cite{Ditor}, \cite{Pelczynski}, that a positive operator
$\L:C(X)\to C(X)$ is continuous iff the mapping $X \ni x
\longmapsto \mu_x \in  \M(X)$ where $\mu_x(a):=\L(a)(x)$, $a
\in C(X)$, is continuous. Furthermore,  if we assume  that
$\mu_x=0$ for $x \in \overline{X\setminus \al(\Delta)}$, then the
mapping $X \ni x \longmapsto \mu_x \in  \M(X)$ is continuous iff
its restriction to $\Int(\al(\Delta))$ is continuous and vanishing
at  infinity.
\\
In view of the above formula \eqref{transfer equation for
commutative1} defines a bounded positive operator $\L$ which  clearly  satisfies  \eqref{b,,2}.
Conversely, if $\L$ is a transfer operator for $(\A,\alpha)$, then 
$\mu_x=0$ for $x \in \overline{X\setminus \al(\Delta)}$. Indeed,   for every $h\in C(X)$ such that $h(\al(\Delta))=1$ we have
$\alpha(h)=\alpha(1)$, and hence
$$
\mu_x(a)=\L(a)(x)=\L(\alpha(1)a)(x)=h(x)
\L(a)(x)=h(x)\mu_x(a), \qquad a \in \A.
$$
This together with Urysohn lemma proves our claim. Let  $x \in
\Int(\al(\Delta))$. Equation \eqref{b,,2} implies that $\int_{X} a\, d\mu_x =\int_{\Delta} a\,d\mu_x$   and  $\int_{\Delta} (a\circ \al)(y) d\mu_x (y)=\L(1)(x)a(x)$, $a\in C(X)$. Using these one gets  $\supp \mu _x \subset \al^{-1}(x)$, and thus $\L$ is of the form \eqref{transfer equation for
commutative1}.
The items i) and ii) follow immediately from Proposition \ref{non-degenerate operators} and Definition \ref{definition of complete transfers}.
\end{Proof}
As we shall see in Corollary \ref{ostateczne wnioski} and Proposition \ref{transfery do cholery} it often happens that uniqueness of non-degenerate
transfer operators is equivalent to completeness, however in general this is not true.
\begin{prz}
Let $X=[0,1]\cup \{2\}$ and $\al|_{[0,1]}=id$, $\al(2)=0$. For a continuous mapping $[0,1] \ni x \longmapsto \mu_x \in  \M(X)$ such that $\mu_x$ is a probability measure supported on $\al^{-1}(x)$ one clearly must have $\mu_x=\delta_x$. Thus the system $(C(X),\alpha)$,  where $\alpha(a)=a\circ \al$, has  a unique non-degenerate transfer operator, even though it does not admit a complete transfer operator.
\end{prz}
It follows from Proposition \ref{transfer equation for commutative stw1}  that   non-degenerate transfer operators for $(\A,\alpha)$ correspond to  norm one left inverses  of  the composition operator $\alpha:C(\al(\Delta))\to C(\Delta)$, that is to the so-called \emph{regular averaging operators} for $\al:\Delta\to \al(\Delta)$, cf. \cite{Ditor}. We now  adapt the results of \cite{Ditor} to obtain criteria for existence of transfer operators. 
  \begin{defn} if $X,Y$ are compact Hausdorff spaces, we denote by $F(X)$  the collection of all closed  subsets of $X$, and say that a  mapping $\Phi: Y\to F(X)$ is \emph{lower semicontinuous (l.s.c)} if for every open $V\subset X$ the set $\{x \in Y: V \cap \Phi(x)\neq \emptyset\}$ is open in $Y$.   \end{defn}
By a \emph{section} of $\al:\Delta \to X$ we mean   a  mapping $\Phi: \al(\Delta)\to F(\Delta)$
 such that  $\Phi(x)\subset \al^{-1}(x)$  for all  $x\in \al(\Delta)$.
\begin{stw}\label{necessary for non-zero}
A necessary condition that the system $(\A,\alpha)$ admit
\begin{itemize}
\item[i)] a non-zero transfer operator is  that $\al:\Delta\to \al(\Delta)$
admits a l.s.c section  not identically equal to
$\emptyset$.
\item[ii)]  a non-degenerate  transfer operator is  that $\al(\Delta)$ is open
in $X$ and  $\al:\Delta\to \al(\Delta)$ admit a l.s.c section   taking
values in $F(\Delta)\setminus \{\emptyset\}$.
\end{itemize}
\end{stw}
\begin{Proof}
It is known, cf. \cite[Cor. 3.2]{Ditor},  that for any continuous operator
$\L:C(X)\to C(X)$  the mapping $X \ni x
\longmapsto \supp \mu_x \in  F(X)$ where $\mu_x(a):=\L(a)(x)$, $a
\in C(X)$, is l.s.c.  Hence the statement follows from Proposition \ref{transfer equation for commutative stw1}.
\end{Proof}
 \label{cantor map page} As a non-trivial  application of the above result we note that the only lower
semicontinuous section for 
the Cantor map  $\al:[0,1]\to[0,1]$ (which graph is Devil's staircase) is $\Phi\equiv\emptyset$ and hence  there are no
non-zero transfer operators for $(C([0,1]),\alpha)$ where $\alpha(a)=a\circ \gamma$.
For metrizable spaces the conditions from Proposition \ref{necessary for non-zero} are not only necessary but also sufficient.
\begin{thm}\label{necessary sufficient conditions}
If $X$ is a metric space, then   the system $(\A,\alpha)$ admit 
\begin{itemize}
\item[i)] a non-zero transfer operator if and only if $\al:\Delta\to \al(\Delta)$
admits a l.s.c section  not identically equal to
$\emptyset$.
\item[ii)]  a non-degenerate  transfer operator  if and only if $\al(\Delta)$ is open
in $X$ and  $\al:\Delta\to \al(\Delta)$ admit a l.s.c section  taking
values in $F(\Delta)\setminus \{\emptyset\}$.
\end{itemize}
\end{thm}
\begin{Proof}
 To show item i) let $\Phi:\al(\Delta)\to F(\Delta)$ be a l.s.c section for $\al$  not identically equal to
$\emptyset$. Then the set $V:=\{x \in \al(\Delta):  \Phi(x)\neq \emptyset\}$ is open and not empty. By Urysohn lemma there exists  a non-zero continuous  mapping $h:X\to [0,1]$ with a support contained in a  closed subset $K\subset V$. By \cite[Thm. 3.4]{Ditor} there exists a regular averaging operator $U:C(\al^{-1}(K))\to C(K)$  for $\al:\al^{-1}(K)\to K$ and one easily sees that the formula 
$$
\L(a)(x):=\begin{cases}
h(x)(Ua|_{\al^{-1}(K)})(x)&,  x \in
K\\
0   &, x\in X\setminus K
\end{cases},\qquad a\in \A=C(X),
$$
defines a non-zero transfer operator for $(\A,\alpha)$.
Item ii) follows  from Proposition \ref{transfer equation for commutative stw1} i) and \cite[Thm. 3.4]{Ditor}. 
\end{Proof}
We end this  section  discussing the  case  of finite-to-one mappings. To this end for every transfer operator $\L$ we  define a function $\rho: \Delta \to [0,\infty )$ by the formula
\begin{equation}\label{rho function}
\rho(x):=\mu_{\al(x)}(\{x\}), \qquad x\in \Delta,
\end{equation}
where $\mu_{x}\in \M(\Delta)$, $x\in \al(\Delta)$, is given by $\mu_{x}(a)=\L(a)(x)$.  In general, this function is not continuous, does not determine operator $\L$ uniquely, and not every continuous function might serve as $\rho$. We have, however,  the following 
 \begin{stw}\label{transfer equation for commutative stw}
 If the  map $\al$ determined by \eqref{endomorphisms} is finite-to-one, then
every transfer operator for $(\A,\alpha)$ is of the form
\begin{equation}\label{transfer equation for commutative}
\L(a)(x)=\begin{cases}
\sum_{y\in \al^{-1}(x)} \rho(y) a(y)&,  x \in \al(\Delta)\\
0   &, x\notin \al(\Delta)
\end{cases}
\qquad a \in \A,
\end{equation}
where $\rho: \Delta \to [0,\infty )$ is a  function given by \eqref{rho function}.
Moreover,
\begin{itemize}
\item[i)] if there exists an open cover $\{V_i\}_{i=1}^n$ of $\Delta$ such that $\al$ restricted to $V_i$, $i=1,...,n$,  is one-to-one, then $\rho$ is necessarily continuous,
\item[ii)] if item $i)$ holds and $\al:\Delta\to X$ is an open map (that is $\al:\Delta\to \al(\Delta)$ is a local homeomorphism and $\al(\Delta)\subset X$ is open), then every continuous function $\rho: \Delta \to [0,\infty )$ defines via  \eqref{transfer equation for commutative}  a transfer operator for $(\A,\alpha)$.
\end{itemize}
\end{stw}
\begin{Proof} As for each $x\in \al(\Delta)$ the support $\supp \mu_x\subset \al^{-1}(x)$ consists of finite number of points, the formula \eqref{transfer equation for commutative} follows from \eqref{transfer equation for commutative1} and \eqref{rho function}.\\
i).  Let
$\{h_i\}_{i=1}^n\subset C(\Delta)$ be a partition of unity subject to
$\{V_i\}_{i=1}^{n}$. Treating $h_i$ as an element of $C(X)$ we   get
$$
\alpha(\L(h_i))(x)=\L(h_i)(\al(x))=\rho(x)h_i(x),\qquad x\in V_i.
$$
Hence for each $i=1,...,n$  the function $\rho(x)h_i(x)$ is continuous on $\Delta$, and thus 
$
\rho = \sum_{i=1}^n \rho 
\cdot h_i\
$ is  continuous as well. \\
ii). Is straightforward since for the local homeomorphism $\al:\Delta\to \al(\Delta)$, 
the number of elements in $\al^{-1}(x)$ for $x\in \al(\Delta)$ is constant on every component of  $\al(\Delta)$ and $\al(\Delta)$ is clopen. 
\end{Proof}
\begin{wn}\label{ostateczne wnioski}
If $\al:\Delta\to \al(\Delta)$ is a local homeomorphism and $\al(\Delta)\subset X$ is open, then \eqref{transfer equation for commutative} establishes a  one-to-one correspondence between 
the  non-degenerate transfer operators for $(\A,\alpha)$ and  
continuous functions $\rho:\Delta\to [0,\infty)$ such that $\sum_{\al(x)=\al(y)} \rho(y)=1$ for all $x\in \Delta$.
In particular, there is a unique non-degenerate transfer operator for $(\A,\alpha)$ iff it is a complete transfer operator. 
\end{wn}
\section{Transfer operators  for systems on $\B(H)$}\label{transfers on B(H)
section}

Let  $H$ be a separable Hilbert space. We recall that due
to \cite[Thm 1]{Tomiyama}  every bounded positive linear map
$T:\B(H)\to \B(H)$ has a  unique decomposition into its normal and
singular part:
$$
T= T_n + T_s
$$
where $T_n$, $T_s$ are bounded positive linear operators, $T_n$ is
normal and $T_s$ is singular ($T_s$ annihilates  the algebra
$\K(H)$ of compact operators). Moreover, taking into account the
explicit form of Tomiyama's decomposition, see \cite{Tomiyama},
one sees that if  $T$ is an endomorphism (resp. a transfer operator,
or a conditional expectation)  $T_n$, $T_s$ are again endomorphisms
(resp. transfer operators, conditional expectations). In
particular, since there is no non-zero representation of the Calkin
algebra on a separable Hilbert space, every endomorphism $\alpha$
of $\B(H)$ is normal, and the classification of transfer
operators for $(\B(H),\alpha)$  splits into classification of all
normal and all singular ones.
 \\
We fix an endomorphism  $\alpha:\B(H)\to \B(H)$. Normality of $\alpha$ implies that there
exists a  number $n=1,2,..., \infty$  called  \emph{an index} for $\alpha$,
and a family of isometries $U_i$, $i=1,...,n$  with orthogonal ranges 
 such that
\begin{equation}\label{powers enodmorphisms}
 \alpha(a)= \sum_{i=1}^{n} U_i a U_i^*, \qquad a \in \B(H),
\end{equation}
where (if $n=\infty$) the sum is weakly convergent, see, for instance, \cite{Laca}.
As transfer operators for $(\B(H),\alpha)$ are completely determined by their
action on the algebra  $\alpha(1)\B(H)\alpha(1)$
 we note that $\alpha(1)\B(H)\alpha(1)$   may be presented as the following
$W^*$-tensor product
$$
\alpha(1)\B(H)\alpha(1)= \B(H_0)\otimes \B(H_n)=\B(H_0\otimes H_n)
$$
where $\dim H_0 =\infty$ and $\dim H_n =n$. To be more precise,
this decomposition is obtained  assuming the following
identification
\begin{equation}\label{tensor vs matrices}
a\otimes b =U^*aU\sum_{i,j=1}^n b(i,j) U_i U_j^*
\end{equation}
where the sum is weakly convergent, $\{b(i,j)\}_{i,j=1}^n$ is a matrix of
$b\in \B(H_n)$ with respect to a fixed orthonormal basis $\{e_i\}_{i=1}^n$
for $H_n$, and $U$ is the isometry $U_{i_0}$ for  an arbitrarily chosen  but
fixed index  $i_0=1,...,n$  (then $H_0=UU^*H$), cf. \cite{Laca},
\cite{Dixmier}.
Under these identifications, the homomorphism $\alpha:\B(H)\to \B(H_0\otimes H_n)\subset \B(H)$
has an especially simple   form: $\alpha(a)= UaU^* \otimes 1$, $a \in \B(H)$, and in particular $\alpha(\B(H))=\B(H_0)\otimes
1.$
\begin{stw}\label{transfery do cholery}
For every transfer operator $\L$ for $(\B(H),\alpha)$ there exists a
positive linear functional $\omega$ on $\B(H_n)$ such that
\begin{equation}\label{tensoring the transfers}
\L(a \otimes b) = U^*aU \omega(b), \qquad a \otimes b\in \B(H_0)\otimes
\B(H_n).
\end{equation}
Moreover,
\begin{itemize}
\item[i)] relation \eqref{tensoring the transfers} establishes an
isometric isomorphism between the cone of  normal transfer
operators $\L$ for $(\B(H),\alpha)$ and the cone of normal
positive linear functionals $\omega$ on $\B(H_n)$. 
\item[ii)]
$\L$ is singular iff $\omega$ is singular (i.e. annihilates
$\K(H_n)$),  in particular, if $n$ is finite then all  transfer
operators are normal.
 \item[iii)] there is a unique
non-degenerate transfer operator $\L$ for $(\B(H),\alpha)$ iff $n=1$, and then  $\alpha(a)=UaU^*$  and $\L(a)=U^*aU$ is a complete
transfer operator.
\end{itemize}
\end{stw}
\begin{Proof} For any $b\in \B(H_n)$, $\L(1\otimes b)$
commutes with every  $a\in\B(H)$:
\begin{align*}
a \L(1\otimes b)  &= \L(\alpha(a)\cdot 1\otimes b)  =\L\big(
(UaU^*\otimes 1) (1\otimes b)\big)=\L(
UaU^*\otimes b)
\\
&=\L\big( (1\otimes b)  (UaU^*\otimes 1) \big)= \L( (1\otimes
b)\alpha(a))  = \L(1\otimes b)a .
\end{align*}
Thus we have $\L(1\otimes b)=\omega(b)\cdot 1$ for certain $\omega(b)\in \C$, and 
 this relation defines a positive linear functional on $\B(H_n)$ for which
\eqref{tensoring the transfers} holds.\\
 i). That normality of $\L$  implies normality of $\omega$   is straightforward. Conversely, as tensors $a\otimes b$ are
linearly strongly dense in $\B(H_0)\otimes \B(H_n)$,   any normal $\omega$ determines uniquely via \eqref{tensoring the transfers} 
  a normal  operator $\L$. Moreover, since $U$ is an isometry, $\L$ is a transfer operator for $\alpha(a)=UaU^*\otimes 1$ and $\|\L\|=\|\omega\|$.
\\
Item ii)  follows from the equality $\K(H)=\K(H_0)\otimes
\K(H_n)$ and formula \eqref{tensoring the transfers}. 
Item iii)  follows from item i).
\end{Proof}
Combining the above statement with Proposition \ref{conditional expectations} one gets 
\begin{stw}[\cite{Sze-Kai}, Prop. 2.4]\label{wniosek jakis tam}
There is a one-to-one correspondence between the normal
conditional expectations $E:\B(H_0\otimes H_n)\to \B(H_0)\otimes
1$ and normal states $\omega$ on $\B(H_n)$,  established by
the relation
\begin{equation}\label{equation for conditional expectations}
E(a\otimes b)=  a \otimes \omega(b)1, \qquad a\otimes b \in \B(H_0)\otimes
\B(H_n).
\end{equation}
Moreover, if $n$ is finite every conditional   expectation $E:\B(H_0\otimes
H_n)\to \B(H_0)\otimes 1$ is normal.
\end{stw}
Since every normal positive functional $\omega$ on $\B(H_n)$ is
of the form  $\omega(a)=\textrm{Tr}\,(\rho a)$ where $\rho\in
\K(H_n)$ is a positive trace class operator, one may deduce the
following explicit form of normal transfer operators.
  \begin{thm}\label{endmorph of B(H)1}
  An operator $\L: \B(H)\to \B(H)$ is a normal transfer operator for
$(\B(H),\alpha)$ if and only if it is of the form
 \begin{equation}\label{transfer on B(H)1}
  \L(a)=\sum_{i,j=1}^{n}\rho(i,j) U_i^* a U_j, \qquad a \in \B(H)
 \end{equation}
 where $\{U_i\}_{i=1}^{n}$ is a family of isometries satisfying \eqref{powers enodmorphisms} and
  $\{\rho(i,j)\}_{i,j=1}^n$ is a matrix of a positive trace class operator
$\rho\in \K(H_n)$, that is
$$
 \sum_{i,j}\rho(i,j)z_i \overline{z}_j \geq 0 \,\,\,
  \textrm{ for all  }z_i\in \C, \,\,i=1,...,n,\quad \textrm{and}\quad
\textrm{Tr}\, \rho=\sum_{i=1}^n\rho(i,i) < \infty.
$$
   Moreover, $\|\L\|=\textrm{Tr} \, \rho$  and $\L$ is
non-degenerate if and only if $\textrm{Tr} \, \rho=1$.
 \end{thm}
 \begin{Proof}
 It suffices to check  that \eqref{transfer on B(H)1} and \eqref{tensoring the
transfers} agree on simple tensors $a\otimes b \in\B(H)\otimes \B(H_n)$ where
$\omega(a)=\textrm{Tr}\,(\rho a)$. In view of \eqref{tensor vs matrices} it is
straightforward.
\end{Proof}
\begin{wn}
If $n<\infty$, then  every transfer operator for $(\B(H),\alpha)$ is of the form \eqref{transfer on B(H)1} for a positive definite matrix $\{\rho(i,j)\}_{i,j=1}^n$.
\end{wn} 
In order to represent not necessarily normal  transfer operators (in the case  $n=\infty$)
  we shall replace the sum in \eqref{transfer on B(H)1} with a  weak integral. We denote by $\M_{f.a}(\N)$ the set of all positive finite finitely
additive measures on $\N$, and recall that  the Lebesgue integral
$\int_\N f(i)\, d\mu(i)$, $f\in \ell^\infty$, $\mu\in \M_{f.a}(\N)$ allows us to treat elements of $\M_{f.a}(\N)$  as positive functionals  on the algebra $\ell^\infty$ of
bounded sequences on $\N$. 
Moreover, identifying $\ell^\infty$ with the corresponding atomic masa of $\B(H_n)$ one may  use a form of "integration" (usually expressed via  limits along ultrafilters)
to extend 
these functionals from $\ell^\infty$  to $\B(H_n)$, cf.
\cite{Kadison-Singer}, \cite{Anderson}. We shall adopt this method
to transfer operators with the help of the following simple consequence of Riesz lemma.
\begin{lem}
For every  $\mu \in \M_{f.a.}(\N)$  and every bounded function $\phi:\N \to
\B(H)$ there exists a unique bounded operator denoted by $\int_\N\phi(i)\,d\mu(i)$
satisfying
$$
\langle \int_{\N} \phi(i) \, d\mu(i)\, x, y\rangle = \int_{\N} \langle
\phi(i)x, y \rangle \, d\mu(i),\qquad \textrm{ for all }\,\, x,y \in H.
$$
Moreover, $\|\int_{\N} \phi(i) \, d\mu(i)\|=\sup_{i\in \N}\|\phi(i)\|$.
\end{lem}
 \begin{thm}\label{integrating  transfers thm} Let the index of $\alpha$ be $n=\infty$.
For every  $\mu \in \M_{f.a.}(\N)$  the formula
\begin{equation}\label{diagonalizable transfer operators}
\L(a)=\int_{\N} U_i^*aU_i\, d\mu(i),\qquad a\in \B(H)
\end{equation}
defines a transfer operator for $(\B(H),\alpha)$.
Moreover,  $\|\L\|= \mu(\N)$ and
\begin{itemize}
\item[i)] $\L$ is singular iff  $\mu$ has no atoms.
\item[ii)] the positive linear functional $\omega$ associated with $\L$
via \eqref{tensoring the transfers} is diagonalizable, that is   $\omega(a)=\int_{\N}
\langle ae_i, e_i\rangle d\mu(i)$,  $a\in \B(H_n)$,
($\{e_i\}_{i\in \N}$  is the fixed basis).
\item[iii)] $\L$ is non-degenerate iff $\mu(\N)=1$, and then the
conditional expectation $E=\alpha\circ \L$  is given by
$
E(a)=\sum_{j=1}^\infty\int_\N  U_jU_i^* a U_i U_j^* \, d\mu(i).
$
\end{itemize}
\end{thm}
\begin{Proof}
Plainly,  $\L(a)=\int_{\N} U_i^*aU_i\, d\mu(i)$  is  positive. It satisfies \eqref{b,,2} as the ranges of isometries $\{U_i\}_{i\in \N}$ are orthogonal and hence  for $a,b\in \B(H)$ we have
$$
\L(\alpha(a)b)=\int_{\N} U_i^*\left(\sum_{j=1}^{n} U_j a U_j^*\right)b U_i\, d\mu(i)= \int_{\N}  a U_i^*b U_i\, d\mu(i)=a \L(b).
$$
For  $x,y \in H$, $\|x\|=\|y\|=1$, we
have
 $|\langle \L(a)x, y\rangle|\leq \int_{\N} |\langle
U_i^*aU_i x, y \rangle| \, d\mu(i) \leq \|a\|\mu(\N)$, and as
$\L(1) =1\cdot \mu(\N)$
 we get 
$\|\L\|=\mu(\N)$.
 \\
 i). It will follow from ii) and Proposition \ref{transfery do cholery} ii).
 \\
 ii). By \eqref{tensor vs matrices}, for  $a\in \B(H_n)$ we have $1\otimes
a=\sum_{i,j=1}^\infty a(i,j) U_iU_j^*\in \B(H)$ where $a(i,j)=\langle a
e_j,e_i\rangle$. Thus using orthogonality of the ranges of $\{U_i\}_{i\in \N}$  we get
  \begin{align*}
 1\cdot \omega(a)&= \L(1\otimes a)=\int_{\N}
U_i^*\left(\sum_{k,j=1}^\infty a(k,j) U_kU_j^*\right)U_i\, d\mu(i)=\int_{\N} 1\cdot  a(i,i)\, d\mu(i)
\\
&=1\cdot\int_{\N}   a(i,i)\, d\mu(i)=
1 \cdot \int_{\N} \langle ae_i, e_i\rangle d\mu(i).
\end{align*}
  iii). Follows from the equality $\alpha(\L(1))=\mu(\N)\alpha(1)$. 
 \end{Proof}
 \begin{wn}
 For every diagonalizable singular state $\omega$ on $\B(H_n)$ there exists a singular  conditional
expectation $E:\B(H_0\otimes H_n)\to \B(H_0)\otimes 1$ satisfying \eqref{equation for conditional expectations}.
 \end{wn}
 Formula \eqref{diagonalizable transfer operators} describes transfer operators "`diagonalizable with respect to the basis  $\{e_i\}_{i\in \N}$ in $H_n$. Changing this basis  one gets more transfer operators. Namely, for every
matrix $\{u(i,j)\}_{i,j=1}^\infty$ of a unitary operator $u\in \B(H_n)$; $u(i,j)=\langle u
e_j,e_i\rangle$, $i,j\in \N$,  the formula
\begin{equation}\label{transfery zcalkowane}
\L(a)=\sum_{k,j=1}^\infty \int_{\N} u(k,i)\overline{u(k,j)} U_j^*aU_i\,
d\mu(i),\qquad a\in \B(H)
\end{equation}
 defines a transfer operator for $(\B(H),\alpha)$. In particular, taking into account all
unitary matrices and all $\sigma$-additive measure on $\N$, \eqref{transfery zcalkowane} is equivalent to \eqref{transfer on B(H)1}, that is it describes  all normal
transfer operators. However, it is almost certain  that  \eqref{transfery zcalkowane} does
not describe all  transfer operators. The functionals associated
with transfer operators from Theorem \ref{integrating  transfers
thm} are specific extensions of functionals on the algebra of all
diagonalizable operators $\ell^\infty\subset \B(H_n)$, and if   the
Kadison-Singer conjecture is false such extensions even for pure
states might not be unique. Moreover, in the light of the recent
result of C. Akemann and N. Weaver \cite{Akemann-Weaver} there
exist not diagonalizable pure states, i.e.  states not having the
form appearing in  Theorem \ref{integrating  transfers thm}  ii).
\label{pager ref for KS}
In order to clarify the role of Kadison-Singer conjecture in  such
considerations we rephrase Theorem
\ref{integrating  transfers thm}.
\begin{stw}
Let $n=\infty$. For every state $\omega\in (\ell^\infty)^*$ of the atomic masa
$\ell^\infty\subset \B(H_n)$ there exists a non-degenerate transfer operator for
$(\B(H),\alpha)$ satisfying
\begin{equation}\label{problem eqaution}
\L(a\otimes b)=  U^*aU \omega(b)1, \qquad a \in \B(H_0), \, b \in \ell^\infty
\subset \B(H_n),
\end{equation}
and such that, if  $\omega$ is a pure state on $\ell^\infty$ then the functional associated with $\L$ via \eqref{tensoring the transfers} is a pure extension of $\omega$ onto  $\B(H_n)$.
\end{stw}
 We have  two  relevant  problems:
\\[3pt]
   \emph{ Problem 1}: Does equation \eqref{tensoring the transfers} where
$\omega$ is a  pure state of $\B(H_n)$ whose restriction to
$\ell^\infty$ is pure, determines uniquely a non-degenerate transfer
operator ?
\\[3pt]  \emph{ Problem 2}: Does equation \eqref{problem eqaution} where $\omega$ is
a pure state of $\ell^\infty$ determine uniquely
a non-degenerate transfer operator ?
\\[3pt]
 The positive answer to Problem 2  implies the positive answer to
Problem 1. If the  Kadison-Singer conjecture is true then Problems
1 and 2 are equivalent, and finally if the answer to Problem 1 is
positive then Problem 2 and Kadison-Singer problem are equivalent.


\begin{bibdiv}
\begin{biblist}


\bib{aee}{article}{
 author={Abadie, B.},
 author ={Eilers, S.},
 author ={Exel, R.},
  title =  {Morita equivalence for crossed products by Hilbert $C^*$-bimodules},
  journal ={Trans. Amer. Math. Soc.},
  date ={1998},
  volume ={350},
  number={8},
   pages = {3043-3054}
} 




\bib{Akemann-Weaver}{misc}{ 
author={Akemann, C.},
author={Weaver, N.},
title={Not all pure states on $\B(H)$ are
diagonalizable},
status={arXiv:\penalty0
math.OA/0606168}
}





\bib{Anderson}{article}{  author={Anderson, J.},   title =  {Extreme points in sets of positive linear maps on $\B(H)$},  journal ={J. Funct. Anal.},  date ={1979},  volume ={31},  number={2},   pages = {195--217}}


\bib{Ant-Bakht-Leb}{misc}{ 
author={Antonevich, A. B.},
author={Bakhtin, V. I.},
author={Lebedev, A. V.},
title={Crossed product of $C^*$-algebra by an endomorphism,
  coefficient algebras and transfer operators},
status={preprint},
eprint={http://arxiv.org/abs/math/0502415},
note={arXiv:math.OA/0502415},
}



\bib{abl}{misc}{ 
author={Antonevich, A. B.},
author={Bakhtin, V. I.},
author={Lebedev, A. V.},
title={T-entropy and Variational Principle for the spectral radius of transfer and weighted shift operators},
status={preprint},
eprint={http://arxiv.org/abs/0809.3116},
note={arXiv:0809.3116},
}

T-entropy and Variational Principle for the spectral radius of transfer and weighted shift operators
 
 
\bib{Bakht-Leb}{misc}{ 
author={Bakhtin, V. I.},
author={Lebedev, A. V.},
title={When a  $C^*$-algebra is a
coefficient algebra for a given endomorphism},
status={preprint},
eprint={http://arxiv.org/abs/math/0502414},
note={arXiv:math.OA/0502414},
}



\bib{bms}{article}{
 author={Brown, L.G.},
  author={Mingo, J.},
   author={Shen, N.},
   title =  {Quasi-multipliers and embeddings of Hilbert $C^*$-modules},
  journal ={Canad. J. Math.},
  date ={1994},
  volume ={71},
   pages = {1150--1174}
}


 
\bib{Ditor}{article}{
 author={Ditor, S.Z. },
     title =  {Averaging operators in $C(S)$ and semi
continuous sections of continuous maps},
  journal ={ Trans. Amer. Math. Soc. },
  date ={1973},
  volume ={ 175},
   pages = {195-208}
}

\bib{Dixmier}{book}{author ={ Dixmier, J.},
  title =        {Les algebres d'opérateurs dans l'espace hilbertien},
   publisher = {Gauthier-Villars},
  date={1969},
  address ={Paris}
  }

\bib{exel2}{article}{
 author ={Exel, R.},
  title =  {A new look at the crossed-product of a
$C^*$-algebra by an endomorphism},
  journal ={Ergodic Theory Dynam. Systems},
  date ={2003},
  volume ={23},
   pages = {1733-1750}
} 
 
 
 \bib{el}{article}{
 author ={Exel, R.},
  author ={Lopes, A.},
  title =  {$C^*$-algebras, approximately proper equivalence relations, and thermodynamic formalism},
  journal ={Ergodic Theory Dynam. Systems},
  date ={2004},
  volume ={24},
   pages = {1051 - 1082}
} 
 
 



 \bib{fmr}{article}{
  author ={Fowler, N. J.},
 author={Muhly, P. S.}, 
 author ={Raeburn, I.},
  title =  {Representations of Cuntz-Pimsner algebras},
  journal ={Indiana Univ. Math. J.},
  date ={2003},
  volume ={52},
   pages = {569-605}
} 

















\bib{Kadison-Singer}{article}{
 author={Kadison, R. V.},
  author={Singer, I. M.},
   title =  {Extensions of pure states},
  journal ={Amer. J. Math},
  date ={1959},
  volume ={81},
   pages = {383--400}
}
\bib{katsura1}{misc}{ 
author={Katsura, T.},
title={A construction of $C^*$-algebras from $C^*$-correspondences},
series={ Contemp. Math},
date={2003},
organization={Amer. Math. Soc.},
volume={217},
address ={Providence},
pages={173-182}
}

\bib{kwakwa}{misc}{ 
author={Kwa\'sniewski, B. K.},
title={C*-algebras generalizing both relative Cuntz-Pimsner and Doplicher-Roberts algebras},
status={preprint},
eprint={http://arxiv.org/abs/0906.4382},
note={arXiv:0906.4382}
}

\bib{kwa-leb2}{article}{ 
author={Kwa\'sniewski, B. K.},
author={Lebedev, A. V.},
title={Reversible extensions of irreversible dynamical systems: $C^*$-method},
  journal ={Math. Sbor.},
  date ={2008},
  volume ={199},
  number={11},
   pages =     {45-74}}   

\bib{kwa-leb3}{article}{ 
author={Kwa\'sniewski, B. K.},
author={Lebedev, A. V.},
title={Crossed product of a $C^*$- algebra by a semigroup of endomorphisms generated by partial isometries},
  journal ={Integr. Equ. Oper. Theory},
  date ={2009},
  volume ={63},
  pages ={403-425},}


\bib{kwa-leb1}{misc}{ 
author={Kwa\'sniewski, B. K.},
author={Lebedev, A. V.},
title={Relative Cuntz-Pimsner algebras, partial isometric crossed products and reduction of relations},
status={preprint},
eprint={http://front.math.ucdavis.edu/0704.3811},
note={arXiv:0704.3811},
}



\bib{Laca}{article}{
 author={Laca, M.},
     title =  {Endomorphisms of $\B(H)$ and Cuntz algebras},
  journal ={J. Operator Theory},
  date ={1993},
  volume ={30},
   pages = {85--108}
}






\bib{lance}{book}{author ={Lance, E. C.},
  title =        {Hilbert C-*-modules: a toolkit for operator algebraist},
   publisher = {Cambridge University Press},
  date={1995},
  address ={Cambridge},
  }





\bib{Pelczynski}{article}{  author={ Pe\l{}czy\'nski , A.},      title =  {Linear extensions, linear averagings, and their applications to linear topological classification of spaces of continuous functions},   journal ={Dissertationes. Math.},   date ={1968},   volume ={58},    pages = {92--142} }

\bib{Sze-Kai}{article}{
 author={ Sze-Kai  , T.},
     title =  {Index of faithful normal conditional expectations},
  journal ={Proc. Amer. Math. Soc.},
  date ={1991},
  volume ={111},
   pages = {111--118}
}






\bib{Tomiyama}{article}{
 author={Tomiyama  , J.},
     title =  {On projections of norm one in $W^*$-algebras, III},
  journal ={Tohoku Math. J.},
  date ={1959},
  volume ={11},
   pages = {125--129}
}









 
 
 


















 
\end{biblist}

\end{bibdiv}
\end{document}